# A Variation on Heawood List-Coloring for Graphs on Surfaces

**Joan P. Hutchinson**
Department of Mathematics and Computer Science
Macalester College
St. Paul, MN 55105, USA
hutchinson@macalester.edu

**Abstract**

We prove a variation on Heawood-list-coloring for graphs on surfaces, modeled on Thomassen's planar 5-list-coloring theorem. For $\epsilon > 0$ define the Heawood number to be $H(\epsilon) = \lfloor (7 + \sqrt{24\epsilon + 1})/2 \rfloor$. We prove that, except for $\epsilon = 3$, every graph embedded on a surface of Euler genus $\epsilon > 0$ with a distinguished face $F$ can be list-colored when the vertices of $F$ have $(H(\epsilon) - 2)$-lists and all other vertices have $H(\epsilon)$-lists unless the induced subgraph on the vertices of $F$ contains $K_{H(\epsilon)-1}$.

**Keywords:** list-coloring, graph embeddings on surfaces, Heawood number
**AMS Subject Classification:** 05C10, 05C15

## 1 Introduction

Thomassen's celebrated 5-list-coloring theorem for planar graphs, with its proof from "The Book" [1], establishes a stronger result [14]. It proves that if $G$ is a graph embedded in the plane, if $F$ is a face of the embedding, and if the vertices of $F$ have 3-lists and the remaining vertices have 5-lists, then $G$ is list-colorable. We investigate the extent to which an analogous theorem holds for graphs embedded on nonplanar surfaces. The same result does not hold with 3-lists on one face $F$ and 5-lists on the others, for example, when the induced subgraph on the vertices of $F$ contains a complete graph on four or more vertices. We also show that this result does not hold for *locally planar* graphs, graphs embedded with all noncontractible cycles suitably long, which are 5-list-colorable [3].

On the other hand every graph on a surface of Euler genus $\epsilon > 0$ can be $H(\epsilon)$-list-colored (see [7] for history), as well as $H(\epsilon)$-colored [6], where the

Heawood number is given by $H(\epsilon) = \lfloor (7 + \sqrt{24\epsilon + 1})/2 \rfloor$. We prove the following analogue of Thomassen's theorem.

**Theorem 1.1** Let $\epsilon \geq 1$, $\epsilon \neq 3$. Let $G$ be embedded on a surface of Euler genus $\epsilon$ with a distinguished face $F$. Let the vertices of $F$ have $(H(\epsilon) - 2)$-lists and those of $V(G) \setminus V(F)$ have $H(\epsilon)$-lists. Then $G$ is list-colorable unless $G$ contains $K_{H(\epsilon)-1}$ with all its vertices lying on $F$.

In [2] we proved a similar list-coloring result when there are $(H(\epsilon) - 1)$-lists on suitably far apart faces and $H(\epsilon)$-lists elsewhere.

**Proposition 1.2** For infinitely many values of $\epsilon \geq 1$, there is a surface of Euler genus $\epsilon$ and a graph $G$ such that $G$ embeds on that surface with all vertices on one face $F$, $G$ does not contain $K_{H(\epsilon)-2}$, and $G$ is not $(H(\epsilon) - 3)$-list-colorable.

Thus in this context the bound of $H(\epsilon) - 2$ for lists on the vertices of $F$ in Thm. 1.1 is best possible.

## 2 Chromatic and topological background

Recall that a graph can be *k-colored* if one of $k$ colors can be assigned to each vertex so that adjacent vertices receive different colors. Given a graph with a list of colors $L(v)$ for each vertex $v$, the graph can be *L-list-colored* if an element from $L(v)$ can be selected for each $v$ that gives a coloring of $G$. A graph is *k-list colorable* if whenever $|L(v)| \geq k$ for each $v$, the graph can be $L$-list-colored. All chromatic and topological definitions and basic results can be found in [7, 9].

Our results use the next two important theorems on list-coloring; if $v$ is a vertex of $G$, let $deg(v)$ denote the degree of $v$ in $G$.

**Theorem 2.1** [15, 4] a. If $\Delta$ is the maximum degree of a graph $G$, then $G$ can be $\Delta$-list-colored except possibly when $G$ is an odd cycle or a complete graph.

b. If $G$ is a graph with lists $L$ such that $|L(v)| = deg(v)$ for each vertex $v$, then $G$ can be $L$-list-colored except possibly when $G$ is (only) 1-connected and each block (i.e., each maximal 2-connected subgraph) is an odd cycle or a complete graph.

The use of "possibly" in the theorem means that there are lists which will prevent list-colorability; these preventing lists are characterized [15, 4], but are not needed explicitly here.

Let $L$ give lists for each vertex of a graph $G$. Then $G$ is *L-critical* if $G$ is not $L$-list-colorable, but every proper subgraph of $G$ is $L$-list-colorable. When $L(v) = \{1, 2, ..., k-1\}$ for each vertex, $G$ is said to be *k-critical*.

**Theorem 2.2** [8]. Let $G$ be a graph with $n$ vertices and $e$ edges that does not contain $K_k$, and let $L$ be lists for $G$ with $|L(v)| = k - 1$ for every vertex $v$ where $k \geq 4$. If $G$ is $L$-critical, then $2e \geq (k-1)n + k - 3$.

Let $S_\epsilon$ denote a surface of Euler genus $\epsilon$. We use the following facts for a multigraph $G$ embedded on $S_\epsilon$, which are derived from Euler's formula:

$$n - e + f \geq 2 - \epsilon \text{ and } e \leq 3n + 3(\epsilon - 2),$$

when $G$ has $n$ vertices, $e$ edges, $f$ faces and, for the second inequality, when each face has at least three boundary edges.

From the Heawood chromatic bound $H(\epsilon)$, it follows that the largest complete graph that could possibly embed on $S_\epsilon$ is $K_{H(\epsilon)}$. In addition, the least Euler genus $\epsilon$ for which $K_n$ embeds on $S_\epsilon$ is given by the inverse function $\epsilon = \lceil \frac{(n-3)(n-4)}{6} \rceil$ (except that $K_7$ only embeds on the torus and not on the Klein bottle, the two surfaces of Euler genus 2). From this we see that, given $\epsilon > 0$, $K_{H(\epsilon)}$ is the largest complete graph that embeds on $S_\epsilon$ for

$$\left\lceil \frac{(H(\epsilon)-3)(H(\epsilon)-4)}{6} \right\rceil \leq \epsilon \leq \left\lceil \frac{(H(\epsilon)-2)(H(\epsilon)-3)}{6} \right\rceil - 1.$$

For ease of future computation, we spell out these bounds explicitly.

**Lemma 2.3** [2]. For $\epsilon > 0$, set $i = \lfloor \frac{H(\epsilon)-3}{3} \rfloor$ so that $H(\epsilon) = 3i + 3$, $3i + 4$ or $3i + 5$ with $i \geq 1$. Then $K_{H(\epsilon)}$ is the largest complete graph that embeds on $S_\epsilon$ for the following values of $\epsilon$:

a) If $H(\epsilon) = 3i + 3$, then $\frac{3i^2-i}{2} \leq \epsilon \leq \frac{3i^2+i-2}{2}$.

b) If $H(\epsilon) = 3i + 4$, then $\frac{3i^2+i}{2} \leq \epsilon \leq \frac{3i^2+3i}{2}$.

c) If $H(\epsilon) = 3i + 5$, then $\frac{3i^2+3i+2}{2} \leq \epsilon \leq \frac{3i^2+5i}{2}$.

The one exception is the Klein bottle, with $\epsilon = 2$, $H(2) = 7$, $i = 1$, on which $K_6$ is the largest embedding complete graph.

By Lemma 2.3 when $H(\epsilon) \equiv 0, 2 \pmod{3}$ there are $i$ values of $\epsilon$ for which $K_{H(\epsilon)}$ is the largest complete graph embedding on $S_\epsilon$, and when $H(\epsilon) \equiv 1 \pmod{3}$, there are $i + 1$.

We call the cases when $H(\epsilon) = 3i + 4$, $i \geq 1$, and $\epsilon = \frac{3i^2+3i}{2}$, the Special Cases for surfaces and their largest embedding complete graph since only in these cases it is possible numerically, from Euler's formula, to embed $K_{H(\epsilon)+1} - E$, the complete graph on $H(\epsilon) + 1$ vertices minus one edge $E$. When $i = 1$, Ringel proved [10] that $K_8 - E$ does not embed on $S_3$, but when $i > 1$ and $H(\epsilon) \equiv 1, 4, 10 \pmod{12}$ he proved that $K_{H(\epsilon)+1} - E$ does embed in these Special Cases. He left unresolved the other cases when $H(\epsilon) \equiv 7 \pmod{12}$ [11, p. 88].

## 3 Main result

When $G$ is embedded on $S_\epsilon$, $\epsilon > 0$, with a distinguished face $F$ and contains $K_{H(\epsilon)-1}$ with all its vertices lying on $F$ and each having an $(H(\epsilon) - 2)$-list, then clearly $G$ cannot be list-colored, and we say $G$ contains an *F-bad* $K_{H(\epsilon)-1}$. In the following $G_F$ denotes the induced subgraph on $V(F)$. For graphs $G_1$ and $G_2$, $G_1 + G_2$ denotes the join of $G_1$ with $G_2$.

**Lemma 3.1.** For all $\epsilon > 0$, Theorem 1.1 holds for all graphs with $n$ vertices, $n \le H(\epsilon)$, that embed on $S_\epsilon$.

**Proof.** We assume $G$ does not contain an $F$-bad $K_{H(\epsilon)-1}$. If $n \le H(\epsilon) - 2$, then the graph can be $(H(\epsilon) - 2)$-list-colored. When $n = H(\epsilon) - 1$, if $G = K_{H(\epsilon)-1}$, then by assumption it contains a vertex with an $H(\epsilon)$-list, all others with at least an $(H(\epsilon) - 2)$-list, and so can be list-colored. Otherwise $G$ is a proper subgraph of $K_{H(\epsilon)-1}$ with $\Delta \le H(\epsilon) - 2$ and so can be list-colored by Thm. 2.1.

If $n = H(\epsilon)$, suppose at most $H(\epsilon) - 1$ vertices lie on $F$. We can list-color $G_F$ since its maximum degree is at most $H(\epsilon) - 2$ and $G_F$ does not contain $K_{H(\epsilon)-1}$. In addition, for $G$, $\Delta \le H(\epsilon) - 1$, the vertices of $V(G) \setminus V(F)$ each have an $H(\epsilon)$-list, and thus the list-coloring of $G_F$ extends to $G$. Otherwise all $H(\epsilon)$ vertices lie on $F$. We know that for each vertex $x$ on $F$, $G \setminus \{x\}$ is not the complete graph $K_{H(\epsilon)-1}$ so that there are two vertices $y, z$, distinct from $x$, that are not adjacent. Applying this same fact to $G \setminus \{y\}$ shows that $G$ is a subgraph of $K_{H(\epsilon)}$ minus at least two edges. Suppose two missing edges share a common vertex $x$. Then $G \setminus \{x\}$ can be list-colored since it is not $K_{H(\epsilon)-1}$. In $G$, $x$ is adjacent to at most $H(\epsilon) - 3$ vertices and so the list-coloring extends to $x$. Thus we assume $G$ is a subgraph of $K_{H(\epsilon)} \setminus \{ab, cd\}$, defined to be $K_{H(\epsilon)}$ minus two independent edges, $ab$ and $cd$ with $a, b, c,$ and $d$ distinct vertices. $K_{H(\epsilon)} \setminus \{ab, cd\}$ is the join of the 4-cycle $(a, c, b, d)$ with $K_{H(\epsilon)-4}$. $K_{H(\epsilon)-4}$ can be $(H(\epsilon) - 2)$-list-colored. For $x = a, b, c, d$, reset $L(x)$ to be its original list minus the colors on $K_{H(\epsilon)-4}$ so that each $L(x)$ becomes a list of size at least 2.

Then the 4-cycle $(a, c, b, d)$ can always be list-colored by Thm. 2.1, and $K_{H(\epsilon)} \setminus \{ab, cd\}$ is list-colorable, as is every possible subgraph. □

**Lemma 3.2.** If $G$ embeds on $S_\epsilon$, $\epsilon > 0$, with vertices of degree $H(\epsilon) - 2$ and $H(\epsilon)$ and these have, respectively, lists of size $H(\epsilon) - 2$ and $H(\epsilon)$, if all of the former vertices lie on one face $F$, and if the induced subgraph on $V(F)$ does not contain $K_{H(\epsilon)-1}$, then $G$ can be list-colored.

**Proof.** If $G = K_{H(\epsilon)-1}$, by assumption at least one vertex does not lie on $F$ so that $G$ can be list-colored. Otherwise by Thm. 2.1 $G$ can be list-colored unless it is (only) 1-connected, contains at least two blocks, and possibly has coloring-preventing lists. We consider the block-cutvertex tree of $G$, $\mathrm{BCT}(G)$, which consists of a vertex for each cutvertex and a vertex for each block with a cutvertex adjacent to a block-vertex if and only if the cutvertex lies on the block. Let $P^*$ be a maximal path in $\mathrm{BCT}(G)$ and let $v^*, w^*$ be the two endpoints of $P^*$. Then $v^*$ and $w^*$ correspond to two blocks of $G$, say $G_{v^*}$ and $G_{w^*}$. $G_{v^*}$ has one vertex that is a cutvertex of $G$, say $v'$, and $\deg(v') > \deg(v)$ for every $v \neq v'$ in $G_{v^*}$. Thus $\deg(v') = H(\epsilon)$ and $\deg(v) = H(\epsilon) - 2$ for every $v \neq v'$ in $G_{v^*}$. Thus $G_{v^*}$ and similarly $G_{w^*}$ are the complete graphs $K_{H(\epsilon)-1}$. (Since $H(\epsilon) \geq 6$, these blocks cannot be odd cycles.) Every vertex of $G_{v^*}$, except for $v'$, must lie on $F$, and similarly every vertex of $G_{w^*}$, except for $w'$, its cutvertex, must lie on $F$, and $v' \neq w'$ (since $2(H(\epsilon) - 2) > H(\epsilon) = \deg(v') = \deg(w'))$. Thus $G$ contains two disjoint copies of $K_{H(\epsilon)-1}$, and in fact it must contain at least three disjoint copies of $K_{H(\epsilon)-1}$. Vertex $v'$ is incident with two additional edges besides the $H(\epsilon) - 2$ of $G_{v^*}$. These two edges must lie on an odd cycle, for if instead each was an incident $K_2$, then there would be a longer path than $P^*$ in $\mathrm{BCT}(G)$. Thus the vertex $z^*$ of $\mathrm{BCT}(G)$ corresponding to the odd cycle has degree at least three in $\mathrm{BCT}(G)$, and tracing a longest path from $z^*$, that does not lie on $P^*$, yields another leaf of $\mathrm{BCT}(G)$ that represents another (and disjoint) copy of $K_{H(\epsilon)-1}$, say $G_{x^*}$ with cutvertex $x'$ and with all other vertices of degree $H(\epsilon) - 2$ and lying on $F$.

Then $G$ has at least $3(H(\epsilon) - 1)$ vertices and at least $3(H(\epsilon) - 1)(H(\epsilon) - 2)/2 + 3$ edges, and since $e \leq 3n + 3(\epsilon - 2)$, we have $3(H(\epsilon) - 1)(H(\epsilon) - 2)/2 + 3 \leq 9(H(\epsilon) - 1) + 3(\epsilon - 2)$, which is false for $\epsilon \geq 5$ with $H(\epsilon) \geq 9$.

Recall that the Euler genus of a connected graph equals the sum of the genera of its blocks [13]. Suppose $\epsilon = 1$ with $H(\epsilon) = 6$. Then a graph that contains three copies of $K_5$ has Euler genus at least 3 so that three copies of $K_5$ cannot embed on the projective plane. For $\epsilon = 2$ with $H(\epsilon) = 7$, a graph that contains three copies of $K_6$ has Euler genus at least 3 and so does not embed on the torus or Klein bottle. When $\epsilon = 4$ with $H(\epsilon) = 8$, a graph that contains three copies of $K_7$ has Euler genus at least 6 and so does not embed on the double torus or the sphere with four crosscaps. When $\epsilon = 3$, three copies of $K_6$ can embed on $S_\epsilon$, but by Euler's formula they cannot embed, each with 5 of their 6 vertices on the same face $F$. □

**Proposition 3.3.** Theorem 1.1 holds for all $\epsilon > 0$ provided $\epsilon \neq \frac{3i^2 + 3i}{2}$ for some $i \geq 1$.

**Proof.** The proof is by induction on $n$, the number of vertices of $G$, and holds for $n \leq H(\epsilon)$ by Lemma 3.1. We suppose $n \geq H(\epsilon) + 1$. Denote the number of edges of $G$ by $e$ and the number of faces of the embedding by $f$.

Suppose $G$ contains a vertex $v$ of degree at most $H(\epsilon) - 3$. If $v$ does not lie on $F$, then $G \setminus \{v\}$ does not contain an $F$-bad $K_{H(\epsilon)-1}$, by induction $G \setminus \{v\}$ can be list-colored, and the coloring extends to $v$. If $v$ does lie on $F$, then suppose the face $F$ becomes the face $F'$ in $G \setminus \{v\}$. If the vertices of $F'$ are a subset of the vertices of $F$, then $G \setminus \{v\}$ does not contain an $F'$-bad $K_{H(\epsilon)-1}$, and the result follows again by induction. Otherwise the face $F'$ is extended by at least one vertex $v'$, and $v'$ has an $H(\epsilon)$-list. If the induced subgraph on $V(F')$ contains $K_{H(\epsilon)-1}$, that complete graph must contain $v'$ and so does not contain an $F'$-bad $K_{H(\epsilon)-1}$. Thus we may assume $G \setminus \{v\}$ does not contain an $F'$-bad $K_{H(\epsilon)-1}$ in all cases, it can be list-colored by induction, and the coloring extends to $v$. We conclude that all vertices have degree at least $H(\epsilon) - 2$.

Suppose $G$ contains a vertex $v$ of degree less than $H(\epsilon)$ that does not lie on $F$. Then $G \setminus \{v\}$ cannot contain an $F$-bad $K_{H(\epsilon)-1}$, by induction it can be list-colored, and that coloring extends to $v$. Thus we may assume that all vertices of degree $H(\epsilon) - 2$ and $H(\epsilon) - 1$ lie on $F$.

In addition $F$ cannot contain a vertex of degree $H(\epsilon) - 1$ or greater, for suppose it did. Let $d_i$ denote the number of vertices of $G$ of degree $i$, let $d_{\geq i}$ denote the number of vertices of $G$ of degree at least $i$, and let $d^F_{\geq H(\epsilon)}$ denote the number of vertices of $F$ of degree at least $H(\epsilon)$.

Summing the face sizes we have
$$2e \geq |V(F)| + 3(f - 1) = d_{H(\epsilon)-2} + d_{H(\epsilon)-1} + d^F_{\geq H(\epsilon)} + 3(f - 1).$$

Then by Euler's formula $e \leq 3n + 3(\epsilon - 1) - d_{H(\epsilon)-2} - d_{H(\epsilon)-1} - d^F_{\geq H(\epsilon)}$.

Summing the vertex degrees we have

$2e \geq (H(\epsilon) - 2) d_{H(\epsilon)-2} + (H(\epsilon) - 1) d_{H(\epsilon)-1} + H(\epsilon) (n - d_{H(\epsilon)-2} - d_{H(\epsilon)-1})$
$= H(\epsilon) n - 2 d_{H(\epsilon)-2} - d_{H(\epsilon)-1}$.

Combining the previous two inequalities gives

(∗) $\quad d_{H(\epsilon)-1} + 2 d^F_{\geq H(\epsilon)} \leq 6(\epsilon - 1) - (H(\epsilon) - 6) n$
$\leq 6(\epsilon - 1) - (H(\epsilon) - 6)(H(\epsilon) + 1) \leq 0$.

The bound of 0 follows by evaluating, for $i \geq 1$, $H(\epsilon) = 3i + 3$ and $\epsilon \leq \frac{3i^2 + i - 2}{2}$, $H(\epsilon) = 3i + 4$ and $\epsilon \leq \frac{3i^2 + 3i - 2}{2}$, and $H(\epsilon) = 3i + 5$ and $\epsilon \leq \frac{3i^2 + 5i}{2}$. In the Special Cases with $H(\epsilon) = 3i + 4$ and $\epsilon = \frac{3i^2 + 3i}{2}$, the upper bound of (∗) is 4. Since we are not in a Special Case, we conclude that $d_{H(\epsilon)-1} = d^F_{\geq H(\epsilon)} = 0$.

Thus $G$ has vertices of degree $H(\epsilon) - 2$ on $F$ and others of degree at least $H(\epsilon)$ off $F$. With this information we repeat the calculations above to get $2e \geq d_{H(\epsilon)-2} + 3(f - 1)$, implying that $e \leq 3n + 3(\epsilon - 1) - d_{H(\epsilon)-2}$, and $2e \geq (H(\epsilon) - 2) d_{H(\epsilon)-2} + H(\epsilon) d_{H(\epsilon)} + (H(\epsilon) + 1)(n - d_{H(\epsilon)-2} - d_{H(\epsilon)}) = (H(\epsilon) + 1) n - 3 d_{H(\epsilon)-2} - d_{H(\epsilon)}$.

Combining the two previous inequalities gives

(∗∗) $\quad (H(\epsilon) - 5) n \leq 6(\epsilon - 1) + d_{H(\epsilon)-2} + d_{H(\epsilon)}$.

When $\epsilon = 1$ and $H(\epsilon) = 6$, we have $n \leq d_4 + d_6$, and since $n = d_4 + d_6 + d_{\geq 7}$ we have $n = d_4 + d_6$. By Lemma 3.2 the graphs on the projective plane can be list-colored. When $\epsilon = 2$ and $H(\epsilon) = 7$, we have $2n \leq 6 + d_5 + d_7$, and since $n = d_5 + d_7 + d_{\geq 8}$, we have $n \leq d_5 + d_7 + 2 d_{\geq 8} \leq 6$, a contradiction. (When $\epsilon = 3$, we are in an excluded case.) When $\epsilon = 4$ and $H(\epsilon) = 8$, similarly we get $3n \leq 18 + d_6 + d_8$ so that $2n \leq 2 d_6 + 2 d_8 + 3 d_{\geq 9} \leq 18$ and $n \leq 9 = H(\epsilon) + 1$. Thus $n = 9$, $n = d_6 + d_8$, and so the graphs on the surfaces of Euler genus 4 can be list-colored by Lemma 3.2.

When $\epsilon \geq 5$ and $H(\epsilon) \geq 9$, we have from (∗∗)

$(H(\epsilon) - 5) n = (H(\epsilon) - 5)(d_{H(\epsilon)-2} + d_{H(\epsilon)} + d_{\geq H(\epsilon)+1})$
$\leq 6(\epsilon - 1) + d_{H(\epsilon)-2} + d_{H(\epsilon)}$

so that $\quad (H(\epsilon) - 6)(d_{H(\epsilon)-2} + d_{H(\epsilon)}) + (H(\epsilon) - 5) d_{\geq H(\epsilon)+1} \leq 6(\epsilon - 1)$. If

$n = d_{H(\epsilon)-2} + d_{H(\epsilon)}$, then by Lemma 3.2, these graphs can be list-colored.

Otherwise $d_{\geq H(\epsilon)+1} > 0$. Then

(∗∗∗) $\quad (H(\epsilon) - 6)\, n = (H(\epsilon) - 6)\,(d_{H(\epsilon)-2} + d_{H(\epsilon)} + d_{\geq H(\epsilon)+1})$

$\quad\quad\quad < (H(\epsilon) - 6)\,(d_{H(\epsilon)-2} + d_{H(\epsilon)}) + (H(\epsilon) - 5)\, d_{\geq H(\epsilon)+1} \leq 6\,(\epsilon - 1).$

Thus $H(\epsilon) + 1 \leq n < 6\,(\epsilon - 1)/(H(\epsilon) - 6)$ and

$0 < 6\,(\epsilon - 1) - (H(\epsilon) + 1)\,(H(\epsilon) - 6) \leq 0$

by (∗), giving a contradiction. □

**Corollary 3.4.** In each Special Case, we may assume that all vertices have degree at least $H(\epsilon) - 2$, all vertices not on $F$ have degree at least $H(\epsilon)$, and $d_{H(\epsilon)-1} + 2\, d^F_{\geq H(\epsilon)} \leq 4$ where $d^F_{\geq H(\epsilon)}$ is the number of vertices on $F$ of degree at least $H(\epsilon)$.

Note that the proof of Prop. 3.3 holds for all $\epsilon > 0$ up to line (∗) where the special results of Cor. 3.4 hold for the Special Cases.

## 4 The Special Cases

**Lemma 4.1.** If Theorem 1.1 holds also for $n = H(\epsilon) + 1$, then it holds for all graphs with $n$ vertices that embed on $S_\epsilon$, except possibly when $\epsilon = 3$.

**Proof.** We know Thm. 1.1 holds for the non-Special Cases by Prop. 3.3, and we assume we are in a Special Case with $H(\epsilon) = 3i + 4$, $\epsilon = \frac{3i^2 + 3i}{2}$, $i > 1$. If the theorem holds also for $n = H(\epsilon) + 1$, then that case joins the base cases of Lemma 3.1, and we may assume $n \geq H(\epsilon) + 2$ for the induction step. From the proof of Prop. 3.3 up to (∗) we have for $\epsilon > 0$

(∗∗∗∗) $\quad d_{H(\epsilon)-1} + 2\, d^F_{\geq H(\epsilon)} \leq 6\,(\epsilon - 1) - (H(\epsilon) - 6)\, n$
$\quad\quad\quad\quad\quad \leq 6\,(\epsilon - 1) - (H(\epsilon) - 6)\,(H(\epsilon) + 2) \leq 0$

for $i \geq 2$ in the Special Cases. Thus again we conclude that $d_{H(\epsilon)-1} = d^F_{\geq H(\epsilon)} = 0$.

Then the proof proceeds exactly as in the proof of Prop. 3.3. From (∗∗∗) we obtain $(H(\epsilon) - 6)\, n < 6\,(\epsilon - 1)$ so that $H(\epsilon) + 2 \leq n < 6\,(\epsilon - 1)/(H(\epsilon) - 6)$ and $0 < 6\,(\epsilon - 1) - (H(\epsilon) + 2)\,(H(\epsilon) - 6) \leq 0$ by (∗∗∗∗), giving a contradiction, except when $\epsilon = 3$. □

**Lemma 4.2.** Let $G$ have $k$ vertices each with at least an $i$-list for some $i$, $1 \leq i < k$. If $G$ has at most $i$ vertices of degree $i$ or greater, then it can be list-colored.

**Proof.** Label the vertices in nonincreasing order by degree: $v_1, \ldots, v_i, v_{i+1}, \ldots, v_k$. The first at most $i$ vertices can be list-colored sequentially since when each is colored, it has at most $i - 1$ colored neighbors. After that each remaining vertex has degree at most $i - 1$ and with a $i$-list can always be list-colored. □

In the Special Cases, whether or not $K_{H(\epsilon)+1} - E$ embeds on $S_\epsilon$, we know that a graph on this surface with $n = H(\epsilon) + 1$ vertices is a subgraph of $K_{H(\epsilon)+1} - E$.

**Proposition 4.3.** Theorem 1.1 holds for all graphs with $n = H(\epsilon) + 1$ vertices in all Special Cases when $\epsilon \geq 9$ and $H(\epsilon) \geq 10$.

**Proof.** For $\epsilon \geq 9$ with $H(\epsilon) \geq 10$, we have $H(\epsilon) = 3i + 4$ and $\epsilon = \frac{3i^2 + 3i}{2}$ for $i \geq 2$. Let $G$ be embedded on $S_\epsilon$, have $n = H(\epsilon) + 1$ vertices, and a face $F$ with $|V(F)| = j \leq H(\epsilon) + 1$. Then $G$ is a subgraph of $K_{H(\epsilon)+1} - E$, and $G_F$, the induced subgraph on $V(F)$, does not contain $K_{H(\epsilon)-1}$. By Cor. 3.4 all vertices of $F$ have degree at least $H(\epsilon) - 2$ and all of $V(G) \setminus V(F)$ have degree exactly $H(\epsilon)$, the maximum possible degree of $G$. In addition the vertices of $F$ satisfy $0 \leq d_{H(\epsilon)-1} + 2 d^F_{\geq H(\epsilon)} \leq 4$. Since $n = H(\epsilon) + 1$, when $d^F_{\geq H(\epsilon)} > 0$, $d^F_{\geq H(\epsilon)}$ counts the number of vertices of degree $H(\epsilon)$ that lie on $F$.

Suppose $d_{H(\epsilon)-1} + 2 d^F_{\geq H(\epsilon)} = 0$ so that all vertices of $V(F)$ have degree $H(\epsilon) - 2$. Then by Lemma 3.2, $G$ can be list-colored. Thus we suppose $0 < d_{H(\epsilon)-1} + 2 d^F_{\geq H(\epsilon)} \leq 4$.

Since the vertices of $V(G) \setminus V(F)$ have degree $H(\epsilon)$, we have that $G_F$ is a subgraph of $K_j - E$ and $G$ is a subgraph of the join $G_F + K_{H(\epsilon)+1-j}$. We label the vertices of $G$ $v_1, \ldots, v_j, \ldots, v_{H(\epsilon)+1}$ with the first $j$ lying on $F$ and the remaining outside; we know that some pair of vertices of $F$ is not adjacent. Our goal is to $(H(\epsilon) - 2)$-list-color $G_F$ in a way that extends to all of $G$.

Suppose $v_1$ and $v_2$ are not adjacent and that $c$ is a list-coloring of $G_F$. If $c(v_1) = c(v_2)$, then $c$ will extend to all of $G$ by coloring the vertices $v_{j+1}, \ldots, v_{H(\epsilon)}, v_{H(\epsilon)+1}$ sequentially. Similarly if either $c(v_1)$ or $c(v_2)$ does not lie in $L(v_k)$ for some $k$, $j + 1 \leq k \leq H(\epsilon) + 1$, then we can relabel $v_k$ to become

$v_{H(\epsilon)+1}$, and $c$ extends to all of $G$ by coloring sequentially.

Suppose $j \leq H(\epsilon) - 2$. If $L(v_1) \cap L(v_2) \neq \emptyset$, then we choose a common color for $v_1$ and $v_2$ and this extends to $G_F$ with maximum degree at most $H(\epsilon) - 3$. Otherwise, without loss of generality, $L(v_1) = \{1, 2, ..., H(\epsilon) - 2\}$ and $L(v_2) = \{H(\epsilon) - 1, ..., 2H(\epsilon) - 4\}$. Since $2H(\epsilon) - 4 > H(\epsilon)$, there is a color for $v_1$ or for $v_2$ that does not lie in $L(v_{H(\epsilon)+1})$, and so $G_F$ can be so list-colored with that coloring extending to $G$.

We conclude that $j = H(\epsilon) - 1, H(\epsilon)$, or $H(\epsilon) + 1$, giving us three cases, each divided into subcases according as $d_{H(\epsilon)-1} + 2 d^F_{\geq H(\epsilon)} = 1, 2, 3$ or $4$. Since the proof techniques for the subcases of the first two cases are the same, we pick illustrative subcases. The third case is proved using Thm. 2.2.

Case I. $j = H(\epsilon) - 1$. Thus there are two adjacent vertices lying off $F$, both of degree $H(\epsilon)$; each vertex of $F$ is adjacent to both of these and has degree $H(\epsilon) - 2, H(\epsilon) - 1$, or $H(\epsilon)$ in $G$.

Subcase A. $d^F_{\geq H(\epsilon)} = 0$ and $d_{H(\epsilon)-1} = 1$ or $3$, or $d^F_{\geq H(\epsilon)} = 1$ and $d_{H(\epsilon)-1} = 1$. Suppose the first case holds with $d^F_{\geq H(\epsilon)} = 0$ and $d_{H(\epsilon)-1} = 1$. Then $G_F$ is a subgraph of $K_{H(\epsilon)-1} - E$ with $H(\epsilon) - 2$ vertices of degree $H(\epsilon) - 4$ (in $G_F$, though degree $H(\epsilon) - 2$ in $G$) and one of degree $H(\epsilon) - 3$. Consider the complement $G^c_F$ which has $H(\epsilon) - 2$ vertices of degree 2 and one of degree 1: there is no such graph. There is the same parity contradiction for the other two pairs of values.

In the following subcases we will apply Lemma 4.2 to a subgraph with $k = H(\epsilon) - 3 \geq 7$ vertices, each of which has a list of size at least $i = k - 1$ or $k - 2$.

Subcase B. $d^F_{\geq H(\epsilon)} = 0$ and $d_{H(\epsilon)-1} = 2$. As argued in Subcase A, $G^c_F$ has $H(\epsilon) - 3$ vertices of degree 2 and two of degree 1. Thus $G^c_F$ consists of one path and (possibly) some cycles. Let $v_1, v_2$ be the first two vertices of the path. Since $H(\epsilon) - 2 \leq |L(v_1) \cup L(v_2)| \leq 2H(\epsilon) - 4$, either they can be colored in $G_F$ with a common color or one can be colored with a color not in $L(v_{H(\epsilon)+1})$. Consider the remaining $H(\epsilon) - 3$ vertices of $G_F - \{v_1, v_2\}$ which we must color. Depending on whether the path had two or more than two vertices, respectively, $G^c_F - \{v_1, v_2\}$ consists of either cycles or cycles plus a path, resp. For each vertex in $G_F - \{v_1, v_2\}$ we remove the one or two colors on $v_1$ and $v_2$ from each of their lists (when present), leaving lists of size at least $H(\epsilon) - 4$. Either

$H(\epsilon) - 3$ or $H(\epsilon) - 5$ vertices, resp., in $G_F - \{v_1, v_2\}$ have degree at most $H(\epsilon) - 6$ and at most two have degree at least $H(\epsilon) - 5$. By Lemma 4.2 with $k = H(\epsilon) - 3$ and $i = H(\epsilon) - 5$, the coloring of $v_1$ and $v_2$ extends to $G_F$ and then to all of $G$.

Subcase C. $d^F_{\geq H(\epsilon)} = 0$ and $d_{H(\epsilon)-1} = 4$. In this case $G^c_F$ has $H(\epsilon) - 5$ vertices of degree 2 and four of degree 1, and the proof follows that of Subcase B.

Subcase D. $d^F_{\geq H(\epsilon)} = 1$ or 2, respectively, and $d_{H(\epsilon)-1} = 0$. Then $G^c_F$ has $H(\epsilon) - 2$ (resp., $H(\epsilon) - 3$) vertices of degree 2 and one (resp., two) of degree 0. Thus $G^c_F$ consists of at least one cycle; let $v_1, v_2$ be two successive vertices of one such cycle $C^*$. We color as in Subcase B, and consider $G_F - \{v_1, v_2\}$ with the one or two colors of $v_1, v_2$ removed from the lists of these vertices, leaving lists of size at least $H(\epsilon) - 4$. Then at least $H(\epsilon) - 6$ (resp., $H(\epsilon) - 7$) vertices of $G_F - \{v_1, v_2\}$ have degree at most $H(\epsilon) - 6$, at most two vertices (resp., three vertices) have degree $H(\epsilon) - 4$ (there might be 2, resp. 3, when $C^*$ is a 3-cycle), and at most two vertices (when $C^*$ is larger than a 3-cycle) have degree at least $H(\epsilon) - 5$. By Lemma 4.2 with $k = H(\epsilon) - 3$ and $i = H(\epsilon) - 4$, the coloring of $v_1, v_2$ extends to $G_F$ and so to all of $G$.

Subcase E. $d^F_{\geq H(\epsilon)} = 1$ and $d_{H(\epsilon)-1} = 2$. Then $G^c_F$ has $H(\epsilon) - 4$ vertices of degree 2, two of degree 1, one of degree 0, and consists of a path and (possibly) some cycles. The proof follows those of Subcases B and D.

Case II. $j = H(\epsilon)$. This case and its argument is parallel to that of Case I. There is only one vertex lying off $F$, it has degree $H(\epsilon)$, and it is adjacent to every vertex of $F$. The vertices of $F$ have degree $H(\epsilon) - 2$, $H(\epsilon) - 1$, or $H(\epsilon)$ in $G$. When $d^F_{\geq H(\epsilon)} = 0$ and $d_{H(\epsilon)-1} = 1$ or 3, or $d^F_{\geq H(\epsilon)} = 1$ and $d_{H(\epsilon)-1} = 1$, the same parity contradiction as in Subcase A occurs. Then in the remaining subcases we can apply Lemma 4.2 to a subgraph with $k = H(\epsilon) - 2 \geq 8$ vertices, each of which has a list of size at least $k - 1$ or $k - 2$.

When $d^F_{\geq H(\epsilon)} = 0$ and $d_{H(\epsilon)-1} = 2$, the proof proceeds as in Subcase B; the only difference is that either $H(\epsilon) - 2$ or $H(\epsilon) - 4$ vertices in $G_F - \{v_1, v_2\}$ have degree at most $H(\epsilon) - 5$ and at most two have degree at least $H(\epsilon) - 4$. By Lemma 4.2 with $k = H(\epsilon) - 2$ and $i = H(\epsilon) - 4$, the coloring of $v_1$ and $v_2$ extends to $G_F$ and then to all of $G$. When $d^F_{\geq H(\epsilon)} = 0$ and $d_{H(\epsilon)-1} = 4$, the proof technique is the same.

When $d^F_{\geq H(\epsilon)} = 1$ or 2, respectively, and $d_{H(\epsilon)-1} = 0$, the proof is as in Subcase D. Specifically, at least $H(\epsilon) - 5$ (resp., $H(\epsilon) - 6$) vertices have degree at most $H(\epsilon) - 5$ in $G_F - \{v_1, v_2\}$, at most two vertices (resp., three vertices) have degree $H(\epsilon) - 3$, and at most two vertices (two distinct neighbors of $v_1$, $v_2$ on a cycle) have degree at least $H(\epsilon) - 4$. By Lemma 4.2 with $k = H(\epsilon) - 2$ and $i = H(\epsilon) - 3$, the coloring of $v_1$, $v_2$ extends to $G_F$ and so to all of $G$. When $d^F_{\geq H(\epsilon)} = 1$ and $d_{H(\epsilon)-1} = 2$, the proof follows that of previous cases, as in Case I.

**Case III.** $j = H(\epsilon) + 1$. All vertices of $G$ lie on $F$ and we must show that $G$ can be $(H(\epsilon) - 2)$-list-colored. If all vertices have degree $H(\epsilon) - 2$, then it can be list-colored by Lemma 3.2. Otherwise if $G$ cannot be list-colored, it contains $G'$, an $(H(\epsilon) - 2)$-list-critical subgraph, and we apply Thm. 2.2. $G'$ must have at least $H(\epsilon)$ vertices since $G$ does not contain $K_{H(\epsilon)-1}$, the only list-critical graph on $H(\epsilon) - 1$ or fewer vertices.

First suppose that $G'$ has $n' = H(\epsilon)$ vertices, and we apply Thm. 2.2 with $k = H(\epsilon) - 1$. Then $2e' \geq (H(\epsilon) - 2) H(\epsilon) + H(\epsilon) - 4$. The graph $G$ with $H(\epsilon) + 1$ vertices has, in addition, a vertex of degree at least $H(\epsilon) - 2$. Also $H(\epsilon) - 2$ additional diagonals can be added to the face $F$, triangulating that region. Thus for $G$ plus diagonals with $e^*$ edges

$$2e^* \geq H^2(\epsilon) - H(\epsilon) - 4 + 2(H(\epsilon) - 2) + 2(H(\epsilon) - 2) = H^2(\epsilon) + 3H(\epsilon) - 12.$$

We also have $2e^* \leq 6n + 6(\epsilon - 2) = 6(H(\epsilon) + 1) + 6(\epsilon - 2)$. Thus $H^2(\epsilon) - 3H(\epsilon) - 18 \leq 6(\epsilon - 2)$. With $H(\epsilon) = 3i + 4$ and $\epsilon = \frac{3i^2 + 3i}{2}$, we reach a contradiction.

Similarly suppose that $G'$ has $n' = H(\epsilon) + 1$ vertices and we apply Thm. 2.2 as above. Then $2e' \geq (H(\epsilon) - 2)(H(\epsilon) + 1) + H(\epsilon) - 4$. $H(\epsilon) - 2$ additional diagonals can be added to $F$ so that $G$ plus diagonals with $e^*$ edges satisfies $2e^* \geq H^2(\epsilon) + 2H(\epsilon) - 10$ and again with $2e^* \leq 6(H(\epsilon) + 1) + 6(\epsilon - 2)$ we reach a contradiction since $i \geq 2$. □

Then by Prop. 4.3, Lemma 4.1, and Prop. 3.3, we have Thm. 1.1.

## 5 Examples

The only $k$-critical graph on $k + 2$ vertices is the join $K_{k-3} + C_5$ [5]. As shown in [12] this graph with $k = H(\epsilon) - 1$ can sometimes embed on $S_\epsilon$.

**Proof of Prop. 1.2.** Suppose $K_{H(\epsilon)+1} - E$ embeds on $S_\epsilon$ so that we are in a Special Case with $\epsilon = \frac{3i^2 + 3i}{2}$ and $H(\epsilon) = 3i + 4$ (see the comments following

Lemma 2.3). By [10, 11], $i > 1$, and such embeddings are possible when $H(\epsilon) \equiv 1, 4, 10 \pmod{12}$ and perhaps (or perhaps not) when $H(\epsilon) \equiv 7 \pmod{12}$. As described in [12], we can then obtain an embedding of $K_{H(\epsilon)-4} + C_5$ as follows. Suppose $E = v_1 v_2$. Choose three additional vertices $v_3, v_4, v_5$, and delete the edges $v_2 v_3$, $v_3 v_4$, $v_4 v_5$, $v_5 v_1$. The remaining graph is $K_{H(\epsilon)-4} + C_5$. By Euler's formula this graph has too many edges to be embeddable with all vertices on one face (in which face then $H(\epsilon) - 2$ diagonals could be added). Removing one vertex of $K_{H(\epsilon)-4}$ leaves $K_{H(\epsilon)-5} + C_5$ with all vertices on one face. This is an $(H(\epsilon) - 2)$-critical graph which does not contain $K_{H(\epsilon)-2}$ and cannot be $(H(\epsilon) - 3)$-list-colored. □

We ask if there are graphs as in Prop. 1.2 on all surfaces. The next example shows that for all surfaces, 3-lists on one face can prevent list-coloring even without the presence of complete graphs on four or more vertices.

Let $G$ be a 2-connected outerplanar near-triangulation; that is, $G$ is a *triangulated polygon* which consists of a cycle on all vertices plus diagonals that triangulate the interior of the cycle. These graphs are uniquely 3-colorable. Let the vertices on the exterior cycle be labelled successively $v_1, v_2, ..., v_n$. Then for an orientable surface, select two edges on the exterior cycle, say $v_1 v_2$ and $v_{i+1} v_i$, with either $v_1$ and $v_{i+1}$ or $v_2$ and $v_i$ in different color classes of a 3-coloring. After identifying the edge $v_1 v_2$ with $v_{i+1} v_i$, this graph can be embedded on an orientable surfaces with all vertices lying on one face, and the graph cannot be 3-colored or 3-list-colored. For a nonorientable surface the edge $v_1 v_2$ can be identified with a twist with $v_i v_{i+1}$, the resulting graph can be embedded on a nonorientable surface and is not 3-list-colorable when either $v_1$ and $v_i$ or $v_2$ and $v_{i+1}$ lie in different color classes of a 3-coloring of the triangulated polygon. If $n$ is large and the two identified edges are distant in $G$, then the resulting graph on a surface is locally planar and contains no induced $K_4$.

This example also shows that Thm. 1.1 is best possible for the projective plane in terms of the list-size for vertices on the distinguished face $F$. We ask if Thm. 1.1 is similarly best possible for small $\epsilon$ such as $\epsilon = 2$ and 4. We ask if there is a 4-list, 6-list theorem for the Klein bottle, and if there is a 5-list, 7-list theorem for $\epsilon = 3$.

**Acknowledgements** For many helpful mathematical discussions, the author thanks Alice Dean, with appreciation also to D. Archdeacon, A. V. Kostochka, B. Richter, and S. Wagon for their help.